\documentclass[12pt]{article}
\usepackage[pctex32]{graphics}
\usepackage{amsthm,amsmath,amssymb,amscd,verbatim,epsfig}
\usepackage{amssymb}
\usepackage{graphicx}
\usepackage{epsfig}
\usepackage{graphics,pifont}
\newcommand{\beq}{\begin{equation}}
\newcommand{\eeq}{\end{equation}}

\date{}

\newcommand{\f}{\frac}

\newcommand{\lra}{\longrightarrow}

\begin{document}

\title{Obtaining\ New\ Dividing\ Formulas\ $n|Q(n)$\ From\ the\ Known\ Ones}
\author{Bau-Sen Du \\ [.5cm]
Institute of Mathematics \\
Academia Sinica \\
Taipei 11529, Taiwan \\
dubs@math.sinica.edu.tw \\
(Fibonacci Quarterly 38(2000), 217-222) \\}
\maketitle
\begin{abstract}
In this note, we present a few methods (Theorems 1, 2, and 3) from discrete dynamical systems theory of obtaining new functions $Q(n)$ from the known ones so that the dividing formulas $n|Q(n)$ hold.
\end{abstract}


\section{INTRODUCTION}
In {\bf{\cite{lin}}}, Lin introduced a well-known result (i.e. Theorem 3.1, see also {\bf{\cite{br}}}) from discrete dynamical systems theory (which he called iterated maps) concerning the number of period-$n$ points.  As applications, Lin computed the number $N(n)$ of period-$n$  points of the maps $B(\mu, x)$ for some suitably chosen $\mu$ and obtained some interesting dividing formulas $n | N(n)$ (i.e. formulas (4.23) in {\bf{\cite{lin}}}) which was already obtained in {\bf[6}, Theorem 3{\bf ]} from different maps.  As mentioned in {\bf{\cite{lin}}}, each iterated map contributes an $N(n)$, and hence, {\it in principle}, infinitely many $N(n)$ can be obtained.  However, in practice, to actually compute $N(n)$ is not so easy as was demonstrated in {\bf{\cite{lin}}}.  Lin did not mention how to compute explicit formulas for $N(n)$ other than the one for the special maps $B(\mu, x)$ where the method he used does not seem to apply to other maps easily.  In this note, we want to point out that a simple systematic way of constructing functions $Q(n)$ such that $n|Q(n)$ has already been introduced in {\bf{\cite{du1,du2,du3,du4,du5,du6}}} (see also {\bf{\cite{hao}}}) for a large class of {\it {continuous}} maps from a compact interval into itself and examples of various $Q(n)$ can also be found in {\bf{\cite{du3,du4,du5,du6}}}.  Furthermore, we want to present a few methods (Theorems 1, 2, and 3) from discrete dynamical systems theory of obtaining new functions $Q(n)$ from the known ones so that many more $Q(n)$ can be constructed (see, for example, Theorem 4).  Finally, in {\bf{\cite{lin}}}, Lin only considered the numbers of period-$n$ points for iterated maps.  He did not mention the numbers of {\it {symmetric}} period-$(2n)$ points.  Therefore, we also include such examples in Theorem 5.
                                                                             
\section{SOME DEFINITIONS}

Since our main results are taken from discrete dynamical systems theory, we shall use the notations commonly used there (see also {\bf{\cite{du3,du4,du5,du6}}}).  For completeness, we include the definitions of $\Phi_i(\phi, n), i = 1, 2$ below.  Let $\phi(n)$ be an integer-valued function defined on the set of all positive integers.  If $n=p_1^{k_1}p_2^{k_2} \cdots p_r^{k_r}$, where the $p_i$'s are distinct prime numbers, $r$ and $k_i$'s are positive integers, we let $\Phi_1(\phi, 1)=\phi(1)$ and let $\Phi_1(\phi, n) =$
$$
\phi(n)-\sum_{i=1}^r \phi(\f n{p_i})+\sum_{i_1<i_2} \phi(\f n{p_{i_1}p_{i_2}})
- \sum_{i_1<i_2<i_3} \phi(\f n{p_{i_1}p_{i_2}p_{i_3}}) + \cdots 
+ (-1)^r \phi(\f n{p_1p_2 \cdots p_r}),
$$
\noindent
where the summation $\sum_{i_1<i_2< \cdots < i_j}$ is taken over all integers $i_1, i_2, \cdots, i_j$ with $1 \le i_1 < i_2 <$ $\cdots < i_j \le r$.  If $n = 2^{k_0}p_1^{k_1}p_2^{k_2} \cdots p_r^{k_r}$, where the $p_i$'s are distinct odd prime numbers, and $k_0 \ge 0, r\ge 1$, and the $k_i$'s $\ge 1$ are integers, we let $\Phi_2(\phi, n) =$
$$
\phi(n)-\sum_{i=1}^r \phi(\f n{p_i})+\sum_{i_1<i_2} \phi(\f n{p_{i_1}p_{i_2}})
- \sum_{i_1<i_2<i_3} \phi(\f n{p_{i_1}p_{i_2}p_{i_3}}) + \cdots 
+ (-1)^r \phi(\f n{p_1p_2 \cdots p_r}), 
$$
\noindent
If $n = 2^k$, where $k \ge 0$ is an integer, we let $\Phi_2(\phi,n) = \phi(n) - 1$.  

\section{MAIN RESULTS}

Let $S$ be a nonempty set and let $f$ be a function from $S$ into itself.  {\it In the sequel}, for every positive integer $n$, we let $\phi_f(n)$ denote the number (if finite) of distinct solutions of the equation $f^n(x) = x$ in $S$, where $f^n$  denotes the $n^{th}$ iterate of $f: f^1 = f$ and $f^n = f \circ f^{n-1}$ for $n > 1$.  By standard inclusion-exclusion arguments, it is easy to see that, for each positive integer $n$, $\Phi_1(\phi_f, n)$ is the number of periodic points of $f$ with least period $n$.  On the other hand, if $S$ contains the origin and $g$ is an odd function from $S$ into itself,
we let $\psi_g(n)$ denote the number (if finite) of distinct solutions of the equation $g^n(x) = -x$.  In this case, if $g^n(y) = -y$, then $g^{kn}(y) = (g^n)^k(y) = -y$ for every odd integer $k \ge 1$ and $g^{mn}(y) = (g^n)^m(y) = y$ for every even integer $m \ge 1$.  So, it is again easy to see, by the same inclusion-exclusion arguments, that $\Phi_2(\psi_g,n)$ 
is the number of symmetric periodic points (that is, periodic points whose orbits are symmetric with respect to the origin) of $g$ with least period $2n$.  Consequently, we have $\Phi_1(\phi_f, n) \equiv 0$ (mod $n$) and $\Phi_2(\psi_g, n) \equiv 0$ (mod $2n$) for all positive integers $n$.  Therefore, by letting $Q(n) = \Phi_1(\phi_f, n)$ or $Q(n) = \Phi_2(\psi_g,n)$, we obtain that $n|Q(n)$ for all positive integers $n$.  In the following, we shall present a few methods (Theorems 1, 2, and 3) from discrete dynamical systems theory of obtaining new functions $Q(n)$ from the known ones so that many more $Q(n)$ 
can be constructed.

Since $\Phi_1(\phi, n)$ is linear in $\phi$ (note that $\Phi_2(\psi, n)$ is not linear in $\psi$ because of its definition  
on $n = 2^k$), we easily obtain the following result: 

\noindent
{\bf Theorem 1.}
{\it Let $\phi_i, i = 1, 2$, be integer-valued functions defined on the set of all positive integers.  If, for all positive integers $n$, $\Phi_1(\phi_1, n) \equiv 0$ (mod $n$) and $\Phi_1(\phi_2, n) \equiv 0$ (mod $n$), then, for any fixed integers $k$ and $m$, $\Phi_1(k\phi_1 + m\phi_2, n) = k\Phi_1(\phi_1,n) + m\Phi_1(\phi_2, n) \equiv 0$ (mod $n$) for all positive integers $n$.}

Let $f$ and $f_i, 1 \le i \le j$, be functions from $S$ into itself and let $(\Pi_{i=1}^j \phi_{f_i})(n) = \Pi_{i=1}^j 
\phi_{f_i}(n)$ for all positive integers $n$.  If $h$ is a function from $S$ into itself defined by $h(x) = f^k(x)$, then, since $h^n(y) = y$ if and only if $f^{kn}(y) = y$, we obtain that $\phi_h(n) = \phi_f(kn)$.  On the other hand, if $H$ is 
a function from the Cartesian product set $S^j$ into itself defined by $H(x_1, x_2,\cdots, x_j) = (f_1(x_1), f_2(x_2), \cdots, f_j(x_j))$, then, since $(y_1, y_2, \cdots, y_j) = H^n(y_1, y_2, \cdots, y_j) = (f_1^n(y_1), f_2^n(y_2)$, $\cdots, f_j^n(y_j))$ if and only if $y_i = f_i^n(y_i)$ for all $1 \le i \le j$, we obtain that $\phi_H(n)= (\Pi_{i=1}^j \phi_{f_i})(n)$.  If $S$ contains the origin and all $f$ and $f_i, 1 \le i \le j$, are also odd functions, then so are $h$ (when $k$ is odd) and $H$.  Arguments similar to the above also show that $\psi_H(n) = (\Pi_{i=1}^j \psi_{f_i})(n)$ $= \Pi_{i=1}^j \psi_{f_i}(n)$.  Therefore, we obtain the following results: 

\noindent
{\bf Theorem 2.}
{\it Let $f$ and $f_i, 1 \le i \le j$, be functions from $S$ into itself.  Then the following hold:
\begin{itemize}

\item[(a)]
For any fixed positive integer $k$, let $\phi_k(n) = \phi_f(kn)$.  Then $\Phi_1(\phi_k, n) \equiv 0$ (mod $n$) for all positive integers $n$.

\item[(b)]
$\Phi_1(\Pi_{i=1}^j \phi_{f_i}, n) \equiv 0$ (mod $n$) for all positive integers $n$.
\end{itemize}}

\noindent
{\bf Theorem 3.}
{\it Assume that the set $S$ contains the origin and let $g$ and $g_i, 1 \le i \le j$, be odd functions from $S$ into itself.  Then the following hold:
\begin{itemize}

\item[(a)]
For any fixed odd integer $k > 0$, let $\psi_k(n) = \psi_g(kn)$.  Then $\Phi_2(\psi_k, n) \equiv 0$ (mod $2n$) for all positive integers $n$.

\item[(b)]
$\Phi_2(\Pi_{i=1}^j \psi_{g_i}, n) \equiv 0$ (mod $2n$) for all positive integers $n$.
\end{itemize}}

\noindent
{\bf Remark.}
Note that, in Theorem 1, we only require $\phi_i$ to satisfy $\Phi_1(\phi_i, n) \equiv 0$ (mod $n$) while in Theorems 2 $\&$ 
3 we require them to be the numbers of (symmetric respectively) periodic points of all periods for some (odd respectively) maps.  It would be interesting to know if these stronger requirements in Theorems 2 $\&$ 3 can be loosened. 

\section{SOME EXAMPLES}

In {\bf{\cite{du5}}}, we show that, for any fixed integer $j \ge 2$, if $\phi_j(n) = 2^n - 1$ for $1 \le n \le j$ and $\phi_j(n) = \sum_{i=1}^j \phi_j(n-i)$ for $j < n$, then $\phi_j$ satisfies the congruence identities $\Phi_1(\phi_j, n) \equiv 0$ (mod $n$) for all positive integers $n$.  Since the constant functions also satisfy the same congruence identities, it follows from Theorem 1 that, for any fixed integers $j, k$, and $m$ with $j \ge 2$, if $\phi_{j,k,m}(n) = m\phi_j(n) + k$ for all positive integers $n$, then $\Phi_1(\phi_{j,k,m}, n) \equiv 0$ (mod $n$) for all positive integers $n$.  Since it is easy to see that $\phi_{j,k, m}$ also satisfies the recursive formula : $\phi_{j,k,m}(n) = m(2^n-1)+k$ for $1 \le n \le j$ and $\phi_{j,k,m}(n) =\biggl(\sum_{i=1}^j \phi_{j,k,m}(n-i)\biggr) - (j-1)k$ for $j < n$, we have the following 
result: 

\noindent
{\bf Theorem 4.}
{\it For any fixed integers $j$, $k$, and $m$ with $j \ge 2$, let

$$
\phi_{j,k,m}(n) = 
\begin{cases}
m(2^n - 1) + k, & \text{for} \,\,\, 1 \le n \le j, \cr
\biggl(\sum_{i=1}^j \phi_{j,k,m}(n-i)\biggr) - (j-1)k, & \text{for} \,\,\, j < n. \cr
\end{cases}
$$

\noindent
Then $\Phi_1(\phi_{j,k,m}, n) \equiv 0$ (mod $n$) for all positive integers $n$.}

The following is an example of $\Phi_2(\psi, n) \equiv 0$ (mod $2n$).  For other examples, see {\bf{\cite{du4,du6}}}.  By Theorem 3 above, many more examples can be easily generated from these known ones.

\noindent
{\bf Theorem 5.}
{\it Let $j \ge 2$ be a fixed integer and let $g_j(x)$ be the continuous map from $[-j, j]$ onto itself defined by
$$
g_j(x) =
\begin{cases}
x+1, & \text{for} \,\,\, -j \le x \le -2,\cr
j, & \text{for} \,\,\, x = -1,\cr
-j, & \text{for} \,\,\, x = 1, \cr
x - 1, & \text{for} \,\,\, 2 \le x \le j, \cr
\text{linear}, & \text{on each of the intervals} \,\,\, [-2, -1], [-1, 1], [1, 2]. \cr
\end{cases}
$$

\noindent
We let $\phi_j(n)$ be defined by
$$
\phi_j(n) =
\begin{cases}
3^n - 2, & \text{for} \,\,\, 1 \le n \le j, \cr
3^n- 2 - 4n \cdot 3^{n-j-1}, & \text{for} \,\,\, j+1 \le n \le 2j - 1, \cr
\sum_{i=1}^j (2i-1)\phi_j(n-i) + \sum_{i=j+1}^{2j-1} 
(4j-2i-1)\phi_j(n-i), & \text{for} \,\,\, 2j \le n. \cr
\end{cases} 
$$

\noindent
We also let $\psi_j(n)$ be defined by
$$
\psi_j(n) = 
\begin{cases}
3^n, & \text{for} \,\,\, 1 \le n \le j-1,\cr
3^j - 2j, & \text{for} \,\,\, n = j,\cr
3^n - 4n \cdot 3^{n-j-1}, & \text{for} \,\,\, j+1 \le n \le 2j-1,\cr
\sum_{i=1}^j (2i-1) \psi_j(n-i) + \sum_{i=j+1}^{2j-1} (4j-2i-1)
\psi_j(n-i), & \text{for} \,\,\, 2j \le n. \cr
\end{cases}
$$

\noindent
Then, for any integer $j \ge 2$, the following hold: 
\begin{itemize}

\item[(a)]
For any positive integer $n$, $\phi_j(n)$ is the number of distinct solutions of the equation $g_j^n(x) = x$ in $[-j, j]$.  Consequently, $\Phi_1(\phi_j, n) \equiv 0$ (mod $n$) for all positive integers $n$.

\item[(b)]
For any positive integer $n$, $\psi_j(n)$ is the number of distinct solutions of the equation $g_j^n(x) = -x$ in $[-j, j]$.  Consequently, $\Phi_2(\psi_j, n) \! \equiv \! 0$ (mod $2n$) for all positive integers $n$.
\end{itemize}}

\noindent
{\bf Remark.}
Numerical computations suggest that the functions $\psi_j(n)$ in Theorem 5 also satisfy $\Phi_1(\psi_j, n) \equiv 0$ (mod  $n$) for all positive integers $n$.  However, we are unable to verify this.

\section{OUTLINE OF THE PROOF OF THEOREM 5}

The proof of Theorem 5 is based on the method of symbolic representations which is simple and easy to use.  For a description of this method, we refer to, say, {\bf [6}, Section 2{\bf ]}.  Here we only give an outline of the proof.  We 
shall also use the terminology introduced there.  In the following, we shall assume that $j > 2$.  The case $j = 2$ can be proved similarly.

\noindent
{\bf Lemma 6.}
{\it Under $g_j$, we have
$$
\begin{cases}
(-j)1 & \lra (-(j-1))(-(j-2)) \cdots (-3)(-2)(-1)j(-j), \cr
1(-j) & \lra (-j)j(-1)(-2)(-3) \cdots (-(j-2))(-(j-1)), \cr
(i-1)i & \lra i(i+1) \quad  \text{and} \quad  i(i-1) \lra (i+1)i, \quad \text{for} \,\,\, -(j-2) \le i \le -2, \cr
(-2)(-1) & \lra (-1)j \quad  \text{and} \quad  (-1)(-2) \lra j(-1), \cr
(-j)j & \lra (-(j-1))(-(j-2))\cdots (-3)(-2)(-1)j(-j)123 \cdots (j-2)(j-1), \cr
j(-j) & \lra (j-1)(j-2) \cdots 321(-j)j(-1)(-2)(-3) \cdots (-(j-2))(-(j-1)), \cr
12 & \lra (-j)1 \quad \text{and} \quad 21 \lra 1(-j), \cr
i(i+1) & \lra (i-1)i \quad \text{and} \quad (i+1)i \lra i(i-1) \quad \text{for} \,\,\, 2 \le i \le j-2, \cr
j(-1) & \lra (j-1)(j-2) \cdots 321(-j)j, \cr
(-1)j & \lra j(-j)123 \cdots (j-2)(j-1). \cr
\end{cases}
$$}

In the following, when we say the representation for $y = g_j^n(x)$, we mean the representation obtained, following the procedure as described in Section 2 of {\bf [6]}, by applying Lemma 6 to the representation $(-(j-1))(-(j-2))  
\cdots (-3)(-2)(-1)j(-j)123 \cdots (j-2)(j-1)$ for $y = g_j(x)$ successively until we get to the one for $y = g_j^n(x)$.

For every positive integer $n$ and all integers $k, i$ with $-(j-1) \le k \le j-1$ and $-(j-1) \le i \le j-1$, let $a_{n,k,i,j}$ denote the number of $uv$'s and $vu$'s in the representation for $y = g_j^n(x)$ whose corresponding $x$-coordinates are in the interval $[s_k, t_k]$, where
$$
[s_k, t_k] =
\begin{cases}
[k-1, k], & \text{for} \,\,\, -(j-1) \le k \le -1, \cr
[-1, 1], & \text{for} \,\,\, k = 0, \cr
[k, k+1], & \text{for} \,\,\, 1 \le k \le j-1, \cr
\end{cases}
$$

\noindent
and

$$
uv =
\begin{cases}
(-j)1, & \text{for} \,\,\, i = -(j-1),\cr
(i-1)i, & \text{for} \,\,\, -(j-2) \le i \le -1,\cr
(-j)j, & \text{for} \,\,\, i = 0, \cr
i(i+1), & \text{for} \,\,\, 1 \le i \le j-2,\cr
j(-1), & \text{for} \,\,\, i = j-1.\cr
\end{cases}
$$

\noindent
We also define $c_{n,j}$ and $d_{n,j}$ by letting
$$
c_{n,j} = \sum_{k=-(j-1)}^{j-1} a_{n,k,k,j} + \sum_{k=1}^{j-1}
(a_{n,-k,0,j} + a_{n,k,0,j}) + \sum_{k=0}^{j-2} (a_{n,-k,-(j-1),j} +
a_{n,k,j-1,j})
$$
\noindent
and
$$
d_{n,j} = \sum_{k=-(j-1)}^{j-1} a_{n,k,-k,j} + \sum_{k=1}^{j-1}
(a_{n,-k,0,j} + a_{n,k,0,j}) + \sum_{k=0}^{j-2} (a_{n,k,-(j-1),j} +
a_{n,-k,j-1,j}).
$$

It is easy to see that, for every positive integer $n$, $c_{n,j}$ is the number of distinct solutions of the equation $g_j^n(x) = x$ and $d_{n,j}$ is the number of distinct solutions of the equation $g_j^n(x) = -x$.

Now from Lemma 6 above, we find that these sequences $<a_{n,k,i,j}>$ can be computed recursively.

\noindent
{\bf Lemma 7.}
{\it For every positive integer $n$ and all integers $k$ with $-(j-1) \le k \le j-1$, we have

$$
\begin{cases}
a_{n+1,k,-(j-1),j} &= a_{n,k,0,j} + a_{n,k,1,j} + a_{n,k,j-1,j}, \cr
a_{n+1,k,-(j-2),j} &= a_{n,k,0,j} + a_{n,k,-(j-1),j}, \cr
a_{n+1,k,i,j} &= a_{n,k,i-1,j} + a_{n,k,0,j} + a_{n,k,-(j-1),j}, \quad 
\text{$-(j-3) \le i \le -1$}, \cr
a_{n+1,k,0,j} &= a_{n,k,-(j-1),j} + a_{n,k,0,j} + a_{n,k,j-1,j}, \cr
a_{n+1,k,i,j} &= a_{n,k,0,j} + a_{n,k,i+1,j} + a_{n,k,j-1,j}, \qquad \quad 
\text{$1 \le i \le j-3$} ,\cr
a_{n+1,k,j-2,j} &= a_{n,k,0,j} + a_{n,k,j-1,j}, \cr
a_{n+1,k,j-1,j} &= a_{n,k,-(j-1),j} + a_{n,k,-1,j} + a_{n,k,0,j}. \cr
\end{cases}
$$}

The initial values of $a_{n,k,i,j}$ can be found easily as follows:

$$
\begin{cases}
a_{1,k,k+1,j} = 1, & \text{for} \,\,\, -(j-1) \le k \le -2, \cr
a_{1,-1,j-1,j} = 1, \cr
a_{1,0,0,j} = 1, \cr
a_{1,1,-(j-1),j} = 1, \cr
a_{1,k,k-1,j} = 1, & \text{for} \,\,\, 2 \le k \le j-1,  \cr
a_{1,k,i,j} = 0, & \text{elsewhere}. \cr
\end{cases}
$$

Since the initial values of the $a_{1,k,i,j}$'s are known, it follows from Lemma 7, by direct but somewhat tedious computations for $n$ ranging from 1 to $2j$, that we can find explicit expressions (omitted) for the sequences $<a_{n,k,i,j}>$, $-(j-1) \le k \le j-1$, $-(j-1) \le i \le j-1$, $1 \le n \le 2j$, and from there we obtain the following two results:
\begin{itemize}

\item[(a)]
$c_{m,j} = \phi_j(m)$ and $d_{m,j} = \psi_j(m)$ for $1 \le m \le 2j-1$. 

\item[(b)]
$a_{2j,k,i,j} = \sum_{m=1}^j (2m-1)a_{2j-m,k,i,j} + \sum_{m=j+1}^{2j-1} (4j-2m-1)a_{2j-m,k,i,j}$ for all $-(j-1) \le k \le j-1, -(j-1) \le i \le j-1$. 
\end{itemize}

Since, for fixed integers $k$ and $i$ with $-(j-1) \le k \le j-1, -(j-1) \le i \le j-1$, $a_{n,k,i,j}$ is a linear combination of $a_{n-1,k,m,j}, -(j-1) \le m \le j-1$, it follows from part (b) above that $a_{n,k,i,j} = \sum_{m=1}^j (2m-1)a_{n-m,k,i,j} + \sum_{m=j+1}^{2j-1} (4j-2m-1)a_{n-m,k,i,j}$ for all $n \ge 2j$.  Since both $c_{n,j}$ and $d_{n,j}$ are linear combinations of $a_{n,k,i,j}$'s, we obtain that $c_{n,j} = \sum_{m=1}^j (2m-1)c_{n-m,j} + \sum_{m=j+1}^{2j-1} 
(4j-2m-1)c_{n-m,j}$ and $d_{n,j} = \sum_{m=1}^j (2m-1)d_{n-m,j} + \sum_{m=j+1}^{2j-1} (4j-2m-1)d_{n-m,j}$ for all $n \ge 
2j$.  This completes the proof of Theorem 5.

\section{ACKNOWLEDGMENTS}
The author is very indebted to Professor Peter Jau-Shyong Shiue and the anonymous referee for their many valuable suggestions that led to a more desirable presentation of this paper.

\end{document}